\documentclass[12pt,letterpaper]{amsart}

\usepackage{amscd}

\usepackage{amssymb,latexsym,amsmath,amsthm,graphicx}
\numberwithin{equation}{section}

\newtheorem*{remark}{Remark}

\usepackage{xcolor}

\newcommand{\ol}{\overline}

\begin{document}

\title[minimal graphs]{Level Curves of Minimal Graphs} 

\author[Allen Weitsman]{Allen Weitsman}

\address{Email: weitsman@purdue.edu}


\begin{abstract}
We consider minimal graphs $u = u(x,y) > 0$ over  domains $D \subset R^2$ bounded by an unbounded
Jordan arc $\gamma$ on which
$u = 0$. We prove an inequality on the curvature of the level curves of $u$, and prove that  if $D$ is concave, then the sets $u(x,y)>C\ (C>0)$ are  all concave.   A consequence of this is that 
solutions, in the case where $D$ is concave, are also superharmonic.

\end{abstract}

\maketitle

\section{INTRODUCTION} Let $D$ be a  plane domain bounded by an unbounded Jordan arc $\gamma$. 
  In this paper we consider the boundary value problem
for the minimal surface equation
\begin{equation}
\label{eq:bdryvalueprob}
\left\{
\begin{aligned}
  &\text{div} \frac{\nabla u}{\sqrt{1+|\nabla u|^2}}=0 \quad \text{and } u>0 \quad \text{in}\ D\\
        & u=0\quad \text{on}\ \gamma
 \end{aligned}\right.
 \end{equation} 
 We shall study the curvature $\kappa = \pm |d\varphi /ds|$ for level curves $u=C\ \ (C>0)$ where $\varphi$ is the angle of the tangent vector to the curve, and the sign will be taken to be $+$ when the curve bends away from the set where $u>C$.  
 
 {\bf Theorem 1.}  \emph{There exists a constant $K$ depending on $u$ such that, if $u$ as in (\ref{eq:bdryvalueprob}) and $C>0$,
 the curvature $\kappa=\kappa(C)$ of the level curve $u=C$ satisfies the inequality}
 \begin{equation}
 \label{thm1}
 |\kappa |\leq \frac KC.
 \end{equation}
 Further comments regarding the constant $K$ are given in \S 6.
 
 Our next result concerns solutions whose domains are concave.  There is a literature  (see \cite{GLR} and references cited there) regarding the propogation of convexity for level curves of solutions to partial differential equations over convex domains.

However, regarding the possible geometry of $D$ in (\ref{eq:bdryvalueprob}), it follows from a theorem of Nitsche \cite[p.256]{Nitsche} that
$D$ cannot be convex unless $D$ is a halfplane since (\ref{eq:bdryvalueprob}) cannot have nontrivial 
solutions over domains contained in a sector of opening less than $\pi$.  On the other hand, amongst
the examples given in \cite{LW}, there is a continuum of graphs which do have concave domains; specifically 
those given parametrically in the right half plane ${\bf H}$ by 
\begin{equation}
\label{functions}
 z(\zeta) = (\zeta + 1)^{\gamma} - \frac{1}{\gamma(2 - \gamma)}(\bar{\zeta} + 1)^{2 - \gamma}
\quad  (\zeta \in {\bf H},\ \  1<\gamma<2)
\end{equation}
together with the height function $ 2 \Re e\,  \zeta$.   A concave domain $D$ is taken to be one whose complement is an unbounded convex domain.  The boundary of $D$ is then a curve which bends away from the domain. 

 In \S 6 we will verify that the domains for
the graphs of (\ref{functions}) are concave.
In this note we shall prove the following

{\bf Theorem 2. }\emph{ If $u$ is a solution to (\ref{eq:bdryvalueprob}) with $D$ concave and bounded by a $C^2$
curve $\gamma$, then the sets where $u>C$ are concave for each $C>0$.}

This has the curious consequence

{\bf Corollary.} \emph{ If $u$ is as in  Theorem 2 above, then
$u$ is also superharmonic in $D$.}

\section{PRELIMINARIES}

For a solution $u$ to the minimal surface equation over a simply connected domain $D$ we shall slightly
abuse notation by using $u$ to also denote the solution to (\ref{eq:bdryvalueprob}) when given in parametric
form.  
We shall  make use of the parametrization of the surface given by $u$ in isothermal coordinates using 
Weierstrass functions $\left( x(\zeta), y (\zeta ), u (\zeta) \right)$
with $\zeta$ in the right half plane ${\bf H}$.   Our notation will then be
given by 
\begin{equation}
\label{downstairs}
f(\zeta) = x(\zeta) +iy(\zeta)\quad\zeta = \sigma+i\tau\in {\bf H}.
\end{equation}
Then $f(\zeta )$ is univalent and harmonic, and since $D$ is simply connected it can be written in
the form 
\begin{equation}
\label{decomp1}
f(\zeta) = h(\zeta) + \ol{g(\zeta)}\quad\zeta = \sigma+i\tau\in {\bf H}
\end{equation}
where $h(\zeta )$ and $g(\zeta)$ are analytic in ${\bf H}$, 
\begin{equation}
\label{Dilatation}
|h'(\zeta)|>|g'(\zeta)|,
\end{equation}
and
\begin{equation}
\label{decomp2}
u(\zeta )= 2\Re e\, i\int \sqrt{h'(\zeta )g'(\zeta )}\,d\zeta  .
\end{equation}
(cf. \cite[\S 10.2]{Duren}).

Now, $u(\zeta )$ is harmonic and positive in ${\bf H}$ and vanishes on $\partial {\bf H}$.  Thus,
(cf. \cite[p. 151]{Tsuji}), 
\begin{equation}
\label{height}
u(\zeta )= k_0\,\Re e\, \zeta,
\end{equation}
where $k_0$ is a positive constant.  This with
(\ref{decomp2}) gives
 
\begin{equation}
\label{Dilatation1}
 g'(\zeta) = - \frac{k}{h'(\zeta)}\qquad (k=k_0^2/4).
\end{equation}

Then from (\ref{Dilatation}) we have, in particular, that
\begin{equation}\label{decomp3}
|h'(\zeta)|\geq \sqrt k.
\end{equation}

It follows from (\ref{height}) that the level curves of $u$ can be parametrized by $f(\sigma_0+i\tau)$ for
$ - \infty<\tau<\infty$ and fixed values $\sigma_0$.  Then the curvature $\kappa$ corresponding to height
$\sigma_0$ with the sign convention given at the begining for
$$
\varphi =\arctan(y_\tau/x_\tau)
$$
is given by
\begin{equation}
\label{curvature}
	\kappa =\kappa(\sigma_0, \tau )=\frac{d\varphi}{ds}=\frac{1}{(x_\tau^2+y_\tau^2)^{3/2}}(x_\tau y_{\tau\tau}-y_\tau x_{\tau\tau} ).
\end{equation}

To compute (\ref{curvature}) we use (\ref{downstairs}) and (\ref{Dilatation1}) to write
\begin{equation}
\label{partial1}
x_\tau=\frac{\partial}{\partial \tau }\Re e(h+\overline g)=\Re e\, i(h'-k/ h')=-\Im m(h'-k/ h')
=-(|h'|^2+k)\Im m\frac 1{\overline h'}
\end{equation}
\begin{equation}
\label{partial2}
x_{\tau\tau}=-\frac{\partial}{\partial \tau }\Im m(h'-k/ h')
=-\Re e (h''+ kh''/ h'^2)
\end{equation}
\begin{equation}
\label{partial3}
y_\tau=\frac{\partial}{\partial \tau }\Im m(h+\overline g)=\Im m\, i(h'+k/ h')=\Re e(h'+k/ h')
=(|h'|^2+k)\Re e\frac 1{\overline h'}
\end{equation}
\begin{equation}
\label{partial4}
y_{\tau\tau}=\frac{\partial}{\partial \tau }\Re e(h'+k/ h')
=-\Im m (h''- kh''/ h'^2)
\end{equation}
Substituting (\ref{partial1})-(\ref{partial4}) into (\ref{curvature}) we get
$$
\kappa=\frac {|h'|^3}{4(|h'|^2+k)^2}\left (-(\frac 1{\overline h'}-\frac 1{h'})(h''-k\frac{ h''}{ h'^2}
-\overline h''+k\frac{\overline h''}{\overline  h'^2})
+(\frac 1{\overline h'}+\frac 1{h'})(h''+k\frac{ h''}{ h'^2}
+\overline h''+k\frac{\overline h''}{\overline h'^2})\right)
$$
which simplifies down to
\begin{equation}
\label{kappa}
\kappa =\frac{|h'|}{|h'|^2+k}\,\Re e\,\frac{h''}{h'}.
\end{equation}

Summarizing this, we have 

{\bf Lemma 1.} \emph{ With $u$ as in (\ref{eq:bdryvalueprob}) and $k_0$ as in (\ref{height}), then the locus of
$u=C$ is the set $\zeta=\sigma_0+i\tau $, where $\sigma_0=C/k_0$ and $-\infty<\tau<\infty$.  The curvature 
$\kappa$ at each point of this level set satisfies (\ref{kappa}).}

The  proof of Theorem 2 uses the comparison of $\kappa$ in (\ref{kappa}) with the corresponding
curvature $\kappa_1$ of the image of the line $\sigma_0+i\tau\quad (-\infty<\tau<\infty)$ under $h$.
  Since 
$\arg h'= \Im m \,\log h'$, the formula (\ref{curvature}) gives 
\begin{equation}
\label{kappa1}
\kappa_1=\frac 1{|h'|}\Re e \frac{h''}{h'}.
\end{equation}

\section{PROOF OF THEOREM 1}

Since $f$
in (\ref{decomp1}) is a univalent harmonic mapping, we may convert the estimate from \cite[Lemma 1]{AL}
(cf. also ( \cite[p. 153]{Duren})) for a 
univalent harmonic mapping $F=H+\overline G$ in the unit disk ${\bf U}$ to a mapping of the half plane {\bf H}. 

{\bf Lemma 2.}  \emph{ Let $u$ be as in (\ref{eq:bdryvalueprob}) and $f=h+\overline g$ as in
 (\ref{decomp1}).  Then
 $$
 \bigg|\frac{h''(\zeta)}{h'(\zeta)}\bigg|\leq A/\sigma
 $$
 for some absolute constant $A$.}

{\bf Proof of Lemma 2.}   For the univalent harmonic mapping $F=H+\overline G$ of {\bf U}, the estimate of  \cite{AL} is
 
$$
\bigg|\frac{H''(w)}{H'(w))}\bigg|\leq\frac {A_1}{1-|w|},\quad \ \ w\in {\bf U}
$$
for some absolute constant $A_1$.
Now, for $f(\zeta)=h(\zeta)+\overline{g(\zeta})$, let
$$
F(w)=f\left(\frac{1+w}{1-w}\right),\quad w\in {\bf U}.
$$
Then,
$$
h(\zeta)=H\left(\frac{\zeta-1}{\zeta+1}\right),
$$
$$
h'(\zeta)=H'\left(\frac{\zeta-1}{\zeta+1}\right)\frac 2{(\zeta+1)^2},
$$
and
$$
h''(\zeta)=H''\left(\frac{\zeta-1}{\zeta+1}\right)\frac 4{(\zeta+1)^4}
-H'\left(\frac{\zeta-1}{\zeta+1}\right)\frac 4{(\zeta+1)^3}.
$$
Thus,
$$
\bigg|\frac{h''(\zeta)}{h'(\zeta)}\bigg|\leq\frac 2{|\zeta +1|}\left (\frac 1{|\zeta +1|}\frac{A_1}{1-\left|\frac
{\zeta -1}{\zeta +1}\right |}+1\right)
$$

$$
\leq \frac 2{|\zeta +1|}\left (\frac{A_1}{|\zeta+1|-|\zeta -1|}+1\right )\leq \frac 2{|\zeta +1|}\left (\frac{A_2(|\zeta +1|+|\zeta -1|)}{4\sigma}+1\right )
$$

$$
\leq A/\sigma
$$
for some absolute constant $A$. \qed

{\bf Proof of Theorem 1.}  From Lemma 1, Lemma 2, and (\ref{decomp3}) it follows that, on the level set $u=C$,
\begin{equation}
\label{kappabound}
|\kappa|\leq \frac{ A}{\sqrt{k}C}.
\end{equation}
\qed

\section{PROOF OF THEOREM 2}

For convenience, we dismiss the trivial case where $u$ is planar, and hence we may assume that $h'$ is
nonconstant.

From the given hypothesis, it follows that $\gamma$ must have asymptotic angles in both directions as
$z\to\infty$.  By a  rotation we may assume that the asymptotic tangent vectors have directions $\pm \alpha$
for some $0\leq\alpha\leq\pi/2$.

From the concavity of $D$ and the assumption that the asymptotic tangents to $\gamma$ have angles $\pm \alpha$, it follows that $y_\tau \geq 0$ for $\sigma=0$.  Thus, from (\ref{partial3}) it follows that for $\sigma =0$,  $\Re e\, 1/\overline h'\geq 0$, and 
hence $\Re e\, 1/h'\geq 0$.  Since, by (\ref{decomp3}) $1/h'$ is bounded in ${\bf H}$, this means that 
$\Re e\, 1/h'> 0$ thoughout ${\bf H}$.  This in turn gives
\begin{equation}
\label{hprime}
\Re e\, h'(\zeta)> 0\qquad \zeta\in {\bf H}.
\end{equation}

Let $\psi(\tau )=\arg h'(i\tau )$. It follows from (\ref{kappa}) and (\ref{kappa1})  that $0\leq\kappa_1\not\equiv 0$
on $\partial{\bf H}$ so that

\begin{equation}
\label{psideriv}
\frac{d\psi}{d\tau} =\frac{\partial}{\partial\tau}\Im m(\log h')= \Re e \frac{h''}{h'}\geq0\quad\textrm{when}\ 
\tau = 0.
\end{equation}

By (\ref{hprime})
\begin{equation}
\label {psi}
-\pi/2\leq \psi (\tau)\leq\pi /2.
\end{equation}

Now, $-\pi /2<\Im m(\log h') <\pi/2$ in ${\bf H}$, and in particular is a bounded harmonic 
function in ${\bf H}$.  So for $\zeta =\sigma+i\tau\in {\bf H}$,
$$
\Im m \,\log h'(\zeta)=
\frac\sigma\pi\int_{-\infty}^\infty
\frac{\psi(t)dt}{\sigma^2+(t-\tau)^2}.
$$
Then
$$
\Re e\frac{h''(\zeta) }{h'(\zeta)}= \frac{\partial}{\partial \tau}\Im m\,\log h'(\zeta)=
\frac{\partial}{\partial\tau}\left(\frac\sigma\pi\int_{-\infty}^\infty
\frac{\psi(t)dt}{\sigma^2+(t-\tau)^2}\right)
=\frac{2\sigma}{\pi}\int_{-\infty}^\infty\frac{(t-\tau)\psi(t)dt}{(\sigma^2+(t-\tau)^2)^2}.
$$
An integration by parts yields
$$
\Re e\frac{h''}{h'}=\frac{\sigma}{\pi}\left (\frac{-\psi(t)}{\sigma^2+(t-\tau)^2}
\Big |_{-\infty}^\infty+\int_{-\infty}^\infty\frac{\psi '(t)dt}{\sigma^2+(t-\tau)^2}
 \right).
$$
By (\ref{psi}) it follows that the first term on the right vanishes, and by
(\ref{psideriv}) the second term is positive.   Thus $\kappa_1$ in (\ref{kappa1})
and hence $\kappa$ in (\ref{kappa}) are positive in ${\bf H}$.
\qed
\section{PROOF OF THE COROLLARY}

We may write the minimal surface equation for $u$ as
$$
\frac{\Delta u+F}{|\nabla u|^3}=0
$$
where $F=F(u,x,y)=u_y^2u_{xx}+u_x^2u_{yy}-2u_xu_yu_{xy}$.  

Now, for a given function $v(x,y)>0$ the curvature of the level set $v(x,y)=0$
is given by $F(v,x,y)/|\nabla v|^3$  \cite[p. 72]{Gray} which is positive when
the curve bends away from the interior of the domain.   Since Theorem 2 shows
that the level sets $u=c$ which bound the sets $u>c$ each have positive 
curvature, then applying this to $F(u-c,x,y)$ we find that $\Delta u<0$ and hence $u$
is superharmonic in $D$
\qed
\section{Concluding Remarks.}

For the examples (\ref{functions}) of \S 1, 
$$
\Re e\frac{h''}{h'}=\Re e\,\frac{\gamma -1}{\zeta + 1}>0.
$$
for $1<\gamma <2$ so that by (\ref{kappa}) these have concave domains.  

Furthermore, using 
(\ref{kappa}), this shows that Theorem 1 is sharp.  Regarding the constant $K$ in Theorem 1,  the scaling factor
$k$ in (\ref{kappabound}) is consistent with the fact that $\kappa$ would be rescaled by replacing
$u(x,y)$ by $cu(x/c,y/c)$ for $0<c<\infty$.

\bibliographystyle{amsplain}

\end{document}